\documentclass[11pt]{article}

\usepackage{amsmath}
\usepackage{amssymb}
\usepackage{eucal}

\begin{document}

\baselineskip 16pt

\title{On  $n$-maximal subgroups of finite groups}
\author{Vika A. Kovaleva\\
{\small Department of Mathematics,  Francisk Skorina Gomel State University,}\\
{\small Gomel 246019, Belarus}\\
{\small E-mail: vika.kovalyova@rambler.ru}\\ \\
{ Alexander  N. Skiba}\\
{\small Department of Mathematics,  Francisk Skorina Gomel State University,}\\
{\small Gomel 246019, Belarus}\\
{\small E-mail: alexander.skiba49@gmail.com}}

\date{}
\maketitle

\begin{abstract}  
We describe finite soluble groups 
in which every 
$n$-maximal subgroup is  $\cal F$-subnormal
for some saturated formation $\cal F$.

\end{abstract}

 \let\thefootnoteorig\thefootnote
\renewcommand{\thefootnote}{\empty}


\let\thefootnote\thefootnoteorig

\section{Introduction}

Throughout this paper, all groups are finite and $G$ always denotes
a finite group. 
We use $\cal U$, $\cal N$ and ${\cal N}^{r}$
to denote the class of all supersoluble groups, the class of all nilpotent 
groups and the class of soluble groups of nilpotent length at most
$r$ ($r\geq 1$).
The symbol $\mathbb P$ denotes the set of all primes,
$\pi (G)$  denotes  the  set  of  prime  divisors  of 
the order of $G$. 
If $p$ is a prime, then we use ${\cal G}_p $ to denote the class of 
all $p$-groups.

Let $\cal F$ be a class of groups. If $1\in \cal F$, then we write 
$G^{\cal F}$ to denote the intersection of all normal subgroups 
$N$ of $G$ with  $G/N \in \cal F$. The class $\cal F$ 
is said 
to be a {\sl formation} if either $\cal F=\emptyset$ or $1\in 
\cal F$ and every homomorphic image of  $G/G^{\cal F}$
belongs to $\cal F$ for any group $G$.
The formation $\cal F$ is said to be: {\sl saturated}
if 
$G\in \cal F$ whenever $G/\Phi (G)\in \cal F$ for any group $G$;
{\sl hereditary} if $H\in \cal F$ whenever $G\in \cal F$ 
and $H$ is a subgroup of $G$. A group $G$ is called {\sl $\cal F$-critical}    provided
 $G$ does not belong to $\cal F$ but every proper subgroup of $G$ belongs to
$\cal F$.

For any formation function $f:\mathbb P\to \{$group formation$\}$, the 
symbol $LF(f)$ denotes the collection of all groups $G$ such that either 
$G=1$ or $G\ne 1$ and  $G/C_G(H/K)\in f(p)$ for every chief factor 
$H/K$ of $G$ and every $p\in \pi (H/K)$.
It is well-known that
for any non-empty 
saturated formation 
$\cal F$, there is a unique formation function $F$
such that
${\cal F} =LF(F)$ and $F(p)={\cal G}_pF(p)\subseteq \cal F$ 
for all primes $p$, where ${\cal G}_pF(p)$ is the set of all groups $G$ such that 
$G^{F(p)}\in {\cal G}_p$ (see Proposition 3.8 in \cite[Chapter IV]{1}). The
formation function $F$ is called the {\sl canonical local satellite} of $\cal F$.
A chief factor $H/K$ of 
$G$ is called {\sl ${\cal F}$-central} in $G$ provided $G/C_G(H/K)\in
F(p)$ for all primes $p$ dividing $|H/K|$, otherwise it is called
 {\sl $\cal F$-eccentric}.

Fix some ordering $\phi$ of $\mathbb P$. 
The record $p\phi q$  means that $p$  precedes $q$ in $\phi$ and $p\ne q$.
Recall that 
a group $G$ of order 
$p_1^{\alpha _1}p_2^{\alpha _2}\ldots
p_n^{\alpha _n}$
is called {\sl $\phi$-dispersive} whenever 
$p_{1}\phi p_{2}\phi \ldots \phi p_{n}$
and for every $i$ there is 
a normal subgroup of $G$ of order 
$p_1^{\alpha _1}p_2^{\alpha _2}\ldots
p_i^{\alpha _i}$. 
Furthermore, if $\phi$ is such that $p\phi q$  always implies $p>q$, 
then every $\phi$-dispersive group is called  {\sl Ore dispersive}.

By  definition,  every  
formation  is  {\sl $0$-multiply 
saturated}  and  for $n\geq 1$
a formation $\cal F$  is  called  {\sl $n$-multiply saturated} 
if ${\cal F} = LF(f)$, where every non-empty value of the function $f$
is an $(n-1)$-multiply saturated formation (see \cite{Sk} and \cite{ShemSk}).
In fact, almost saturated formations met in mathematical practice are 
$n$-multiply saturated for every natural $n$.
For example, the formations of all soluble  groups, all nilpotent groups, all $p$-soluble groups, 
all $p$-nilpotent groups, all $p$-closed groups, all $p$-decomposable groups, all 
Ore dispersive groups, all metanilpotent groups are $n$-multiply saturated for 
all $n\geq 1$. Nevertheless, the formations of all supersoluble  groups and all 
$p$-supersoluble groups are saturated, but they are not $2$-multiply 
saturated formations.

Recall that a subgroup $H$ of  $G$ is called a {\sl $2$-maximal}  
({\sl second maximal}) subgroup of $G$ whenever $H$ is 
a maximal subgroup of some maximal subgroup $M$  of $G$. 
Similarly we can define {\sl $3$-maximal subgroups}, 
and so on. 

The  interesting  and  substantial  direction  in  finite  group  theory  
consists  in  studying  the  relations 
between the structure of the group and its 
$n$-maximal subgroups.  
One of the earliest publication  
in  this  direction  is  the  article  of  Huppert \cite{3}  
who  established  the  supersolubility  of  a  group  
$G$ whose 
all  second  maximal  subgroups  are  normal.
In  the  same  article  Huppert  proved  that  if  all  $3$-maximal 
subgroups  of  $G$  are  normal  in  $G$,
then  the  commutator  subgroup  $G'$
of  $G$  is    nilpotent   and  the 
chief   rank  of  $G$  is  at  most  $2$.
These  two  results  were  developed  by  many  authors. 
Among  the 
recent  results  on  $n$-maximal  subgroups  
we  can  mention  \cite{4},  where  the  solubility  of  groups 
is  established  in  which  all  $2$-maximal  subgroups  
enjoy  the  cover-avoidance  property,  and  \cite{5, 6, 7},  where 
new  characterizations  of  supersoluble  groups  in  terms  of  
$2$-maximal  subgroups  were  obtained. 
The 
classification  of  nonnilpotent  groups  whose  all  $2$-maximal  subgroups  are 
$TI$-subgroups  appeared  in  \cite{8}. 
Description was  obtained in \cite{9} of groups  whose  every $3$-maximal 
subgroup permutes  with all maximal 
subgroups.     
The  nonnilpotent  groups  are  described  in  \cite{10}  
in  which  every  two  $3$-maximal  subgroups  are 
permutable.      
The  groups  are  described  in  \cite{11}  whose  all  $3$-maximal  subgroups  are  
$S$-quasinormal,  that 
is, permute  with  all  Sylow  subgroups.
Subsequently  this  result  was  strengthened  in  \cite{12}  to  provide 
a description of the groups whose all $3$-maximal subgroups are subnormal. 

Despite  of  all  these  and  many  other  known  results  
about  $n$-maximal  subgroups,  the  fundamental 
work of Mann \cite{13} still retains its value.  
It studied the structure of groups whose $n$-maximal subgroups 
are subnormal.  Mann proved that if all $n$-maximal subgroups of a soluble group  $G$ 
are subnormal and 
$|\pi (G)| \geq n + 1$,
then $G$ is nilpotent; but if  $|\pi (G)| \geq n - 1$, 
then $G$ is  $\phi$-dispersive for some ordering  $\phi$ of 
$\mathbb P$.  Finally, in the case  $|\pi (G)| = n$ 
Mann described $G$ completely.

Let $\cal F$ be a non-empty formation.
Recall that a subgroup $H$ of a group $G$ is said to be 
{\sl $\cal F$-subnormal in $G$} if either $H=G$ or there exists a chain of subgroups 
$H=H_0<H_1<\ldots <H_{n}=G$ such that $H_{i-1}$ is a maximal subgroup of
$H_i$ and $H_{i}/(H_{i-1})_{H_i}\in \cal F$, for $i=1,\ldots
,n$. 

The main goal of this article is to prove the following formation
analogs of Mann's theorems.

{\bf Theorem A.}
{\sl Let $\cal F$ be an $r$-multiply saturated formation
 such that ${\cal N}\subseteq {\cal F}\subseteq
{{\cal N}}^{r + 1}$ for some $r \geq 0$. 
If every $n$-maximal subgroup of a soluble group $G$
is $\cal F$-subnormal in $G$ and $|\pi (G)|\geq n+r+1$, then $G\in \cal F$.}

{\bf Theorem B.}
{\sl Let ${\cal F}=LF(F)$ be a saturated formation such that 
${\cal N}\subseteq {\cal F}\subseteq \cal U$, where
$F$ is the canonical local satellite of $\cal F$. 
Let $G$ be a soluble group with $|\pi (G)|\geq n+1$.
Then all $n$-maximal subgroups of $G$ 
are $\cal F$-subnormal in $G$ if and only if $G$ is a group of one of 
the following types:}

I. {\sl $G\in \cal F$.}

II. {\sl $G=A\rtimes B$, where $A=G^{\cal F}$  and  $B$ are Hall subgroups of $G$, 
while $G$ is Ore dispersive and satisfies the following: }

(1) {\sl $A$ is either of the form $N_1\times \ldots \times N_t$, 
where each $N_i$ is a minimal normal subgroup of $G$, which is a Sylow 
subgroup of $G$,  for $i=1, \ldots , t$, or a Sylow $p$-subgroup of $G$ of exponent $p$
for some prime $p$ and the commutator subgroup, the Frattini subgroup, and 
the center of $A$ coincide, while $A/\Phi (A)$ is an $\cal F$-eccentric chief factor
of $G$;  }

(2) {\sl every 
$n$-maximal subgroup of $G$
belongs to $\cal F$ and induces on the Sylow $p$-subgroup of 
$A$ an automorphism group which is contained in $F(p)$ 
for every prime divisor $p$ of $|A|$.}

In the proof of Theorem B we often use Theorem A and   the following 
useful fact.

{\bf Theorem C.} {\sl Let $\cal F$ be a hereditary saturated formation
such that every $\cal F$-critical group is soluble and it has a normal Sylow $p$-subgroup
$G_p\ne 1$ for some prime $p$.  Then
every $2$-maximal subgroup of $G$ is $\cal F$-subnormal in $G$
if and only if    either    $G\in \cal F$ or $G$ is an $\cal F$-critical group
and $G^{\cal F}$ is a minimal normal subgroup of $G$.}

{\bf Theorem  D.}
{\sl Let $\cal F$ be a saturated formation such that 
${\cal N}\subseteq {\cal F}\subseteq \cal U$.
If every $n$-maximal subgroup of a soluble group $G$ is
$\cal F$-subnormal in $G$ and $|\pi (G)|\geq n$, then $G$ is 
$\phi$-dispersive for some ordering $\phi$  of $\mathbb P$. }

All unexplained notation and terminology are standard. 
The reader is referred to \cite{1} or  \cite{15} if necessary.

\section{Preliminary Results}

Let $\cal F$ be a non-empty formation.
Recall that a maximal subgroup $H$ of  $G$ is said to be 
{\sl $\cal F$-normal} in $G$ if $G/H_G\in \cal F$,
otherwise it is said to be {\sl $\cal F$-abnormal} in $G$.

We use the following results.

{\bf Lemma 2.1.} 
{\sl Let $\cal F$ be a formation and 
 $H$  an $\cal F$-subnormal subgroup  of $G$. }

(1) {\sl If $\cal F$ is hereditary and $K\leq G$, then $H\cap K$
is an $\cal F$-subnormal subgroup in
$K$} \cite[Lemma 6.1.7(2)]{15}.

(2) {\sl If $N$ is a normal subgroup in $G$, then $HN/N$ is an 
$\cal F$-subnormal subgroup in $G/N$ }
\cite[Lemma 6.1.6(3)]{15}.

(3) {\sl If  $K$ is a subgroup of 
$G$ such that $K$ is $\cal F$-subnormal in $H$, then
$K$ is $\cal F$-subnormal in $G$} \cite[Lemma 6.1.6(1)]{15}.

(4) {\sl If $\cal F$ is hereditary and 
$K$ is a subgroup of $G$ such that $G^{\cal F}\leq K$, then $K$
is $\cal F$-subnormal in $G$} \cite[Lemma 6.1.7(1)]{15}. 

The following lemma is evident.

{\bf Lemma 2.2.} {\sl Let $\cal F$ be a hereditary formation.
If $G\in \cal F$, then every subgroup of $G$ 
is $\cal F$-subnormal in $G$. }

{\bf Lemma 2.3.} {\sl Let $\cal F$ be a hereditary saturated formation.
If every $n$-maximal subgroup of  $G$ is 
$\cal F$-subnormal in $G$, then every $(n-1)$-maximal subgroup of $G$
belongs to $\cal F$ and 
every $(n+1)$-maximal subgroup of 
$G$ is $\cal F$-subnormal in $G$. }

{\bf Proof.} We first show that every $(n-1)$-maximal subgroup of $G$
belongs to $\cal F$. Let $H$ be an $(n-1)$-maximal subgroup of $G$
and  $K$ a maximal subgroup of $H$.
Then $K$ is an $n$-maximal subgroup of $G$ and so, by hypothesis,
$K$ is $\cal F$-subnormal in $G$. Hence $K$ 
is $\cal F$-subnormal in $H$ by Lemma 2.1(1). 
Thus all maximal subgroups of $H$  are $\cal F$-normal in 
$H$. Therefore
$H \in \cal F$ since  $\cal F$ is saturated.

Now, let $E$ be an $(n+1)$-maximal subgroup of $G$,
and let $E_1$ and $E_2$ be an $n$-maximal and an
$(n-1)$-maximal subgroup of $G$, respectively, such that
$E\leq E_1\leq E_2$. Then, by the above, $E_2\in \cal F$, so
$E_1\in \cal F$. Hence $E$ is 
$\cal F$-subnormal in $E_1$ by Lemma 2.2. 
By hypothesis, $E_1$ is
$\cal F$-subnormal in $G$. Therefore $E$ is $\cal F$-subnormal in $G$. The lemma is proved.

{\bf Lemma 2.4} (See \cite[Chapter VI, Theorem 24.2]{14}). 
{\sl Let $\cal F$ be a 
saturated formation and  $G$ a soluble group. If
$G^{\cal F}\ne 1$ and every $\cal F$-abnormal maximal subgroup 
of $G$ belongs to $\cal F$, then the following hold:}

(1) {\sl $G^{\cal F}$ is a $p$-group for some prime $p$;}

(2) {\sl $G^{\cal F}/\Phi (G^{\cal F})$ is an $\cal F$-eccentric 
chief factor of $G$;  }

(3) {\sl if $G^{\cal F}$ is a non-abelian group, 
then the center, commutator subgroup, and
Frattini subgroup of $G$ coincide and are of exponent $p$;}

(4) {\sl if $G^{\cal F}$ is abelian, then $G^{\cal F}$ is elementary;}

(5) {\sl if $p>2$, then $G^{\cal F}$ is of exponent $p$; for $p=2$ the exponent of $G^{\cal F}$
is at most $4$;}

(6) {\sl every pair of $\cal F$-abnormal maximal subgroups of $G$ 
are conjugate in $G$.}

{\bf Lemma 2.5} (See \cite[Chapter VI, Theorem 24.5]{14}). {\sl Let
$\cal F$ be a  saturated formation.
Let $G$ be an $\cal F$-critical group
and $G$ has a normal Sylow $p$-subgroup $G_p\ne 1$ for some prime $p$.
Then:}

(1) {\sl $G_p=G^{\cal F}$;}

(2) {\sl $F(G)=G_p\Phi (G)$;}

(3) {\sl $G_{p'}\cap C_G(G_p/\Phi (G_p))=\Phi (G)\cap G_{p'}$, where $G_{p'}$ is some
complement of $G_p$ in $G$.}

{\bf Lemma 2.6} (See \cite[Chapter VI, Theorems 26.3 and 26.5]{14}). {\sl Let $G$ be 
an $\cal U$-critical group.
Then: }

(1) {\sl $G$ is soluble and $|\pi (G)|\leq 3$;}

(2) {\sl if $G$ is not a Schmidt group, then $G$ is Ore dispersive;}

(3) {\sl $G^{\cal U}$  is the unique normal Sylow subgroup of $G$;}

(4) {\sl if $S$ is a complement of $G^{\cal U}$ in $G$, then $S/S\cap \Phi (G)$ 
is either a primary cyclic group or a Miller-Moreno group.}

Recall that the product of all normal subgroups of a group $G$ whose 
$G$-chief factors are $\cal F$-central in $G$ is called 
{\sl $\cal F$-hypercentre of $G$} and denoted
by $Z_{\cal F}(G)$ \cite[p. 389]{1}.

{\bf Lemma 2.7} (See \cite[Lemma 2.14]{17}). {\sl Let
$\cal F$ be a saturated formation and $F$ 
the canonical local satellite
of $\cal F$.
Let $E$ be a normal $p$-subgroup of a group $G$. Then
$E\leq Z_{\cal F}(G)$ if and only if
$G/C_G(E)\in F(p)$.}

The product ${\cal M}{\cal H}$ of the formations $\cal M$
and $\cal H$ is the class of all groups $G$ such that $G^{\cal H}\in \cal M$.

{\bf Lemma 2.8} (See \cite[Corollary 7.14]{ShemSk}). {\sl The
 product of any two $n$-multiply saturated 
formations is an  $n$-multiply saturated 
formation.}

We shall also need the  following evident lemma.

{\bf Lemma 2.9.} {\sl If $G=AB$, then $G=AB^x$ for all $x\in G$.}

Let $\cal F$ be a class of groups and $t$ a natural number with $t\geq 2$. 
Recall that $\cal F$ is called {\sl $\Sigma _t$-closed} if $\cal F$ 
contains all such groups $G$ that $G$ has subgroups
$H_1, \ldots , H_t$ whose indices are pairwise coprime 
and $H_i\in {\cal F}$, for $i=1, \ldots , t$.

{\bf Lemma 2.10} (See \cite[Chapter I, Lemma 4.11]{14}). 
{\sl Every formation of nilpotent groups is $\Sigma _3$-closed.}

If ${\cal F}=LF(f)$ and $f(p)\subseteq {\cal F}$ for all primes $p$, then
$f$ is called an {\sl integrated local satellite} of $\cal F$.
Let $\cal X$ be a set of groups. 
The symbol $l_n$form$\cal X$ denotes the intersection of all
$n$-multiply saturated 
formations $\cal F$ such that ${\cal X}\subseteq \cal F$.  In view of 
\cite[Remak 3.1.7]{15}, $l_n$form$\cal X$ is  an $n$-multiply saturated 
formation.

{\bf Lemma 2.11} (See \cite[Theorem 8.3]{ShemSk}). {\sl Let ${\cal F}$ be an
 $n$-multiply saturated 
formation. Then $\cal F$ has an integrated local satellite $f$ such that 
$f(p)=l_{n-1}$form$(G/O_{p',p}(G)|G\in {\cal F})$ for all primes $p$.}

{\bf Lemma 2.12} (See \cite[Section 1.4]{Sk1}). {\sl Every $r$-multiply saturated 
formation contained in $ {\cal N}^{r+1}$   is hereditary.}

{\bf Lemma 2.13} (See \cite[p. 35]{14}). {\sl For any ordering $\phi$
 of $\mathbb P$ the class of all 
$\phi$-dispersive groups is a saturated formation.}

{\bf Lemma 2.14} (See \cite[Corollary 1.6]{17}). {\sl Let $\cal F$ be 
a saturated formation containing all nilpotent groups and $E$ a 
normal subgroup of  $G$.
If $E/E\cap \Phi (G)\in \cal F$, then  $E\in \cal F$.}

{\bf Lemma 2.15} (See \cite[Theorem 15.10]{14}). {\sl Let 
$\cal F$ be a saturated formation and $G$ a group such that
$G^{\cal F}$ is nilpotent. Let $H$ and $M$ be subgroups of $G$, 
$H\in \cal F$, $H\leq M$ and $HF(G)=G$. If
$H$ is $\cal F$-subnormal in $M$, then $M\in \cal F$.}

\section{Proof of Theorem A}

First we give  two propositions which may be independently interesting   
since they generalize some known results.

{\bf Proposition  3.1.}    {\sl  
 Suppose that $G=A_{1}A_{2}=A_{2}A_{3}=A_{1}A_{3},$  
where  $A_{1}$, $A_{2}$ and $A_{3}$ are  soluble subgroups of $G$.
 If  the  indices
$|G:N_{G}(A'_{1})|$, $|G:N_{G}(A'_{2})|$, $|G:N_{G}(A'_{3})|$ are pairwise  coprime, then
$G$ is soluble. }

{\bf Corollary  3.2.}    {\sl  
 Suppose that $G=A_{1}A_{2}=A_{2}A_{3}=A_{1}A_{3},$  
where  $A_{1}$, $A_{2}$ and $A_{3}$ are soluble subgroups of $G$.
 If  the  indices
$|G:N_{G}(A_{1})|$, $|G:N_{G}(A_{2})|$, $|G:N_{G}(A_{3})|$ are pairwise  coprime, then
$G$ is soluble. }

{\bf Corollary 3.3} (H. Wielandt).  {\sl  If $G$ has
three soluble subgroups $A_{1}$, $A_{2}$ and $A_{3}$ whose indices
$|G:A_{1}|$, $|G:A_{2}|$, $|G:A_{3}|$ are pairwise  coprime, then
$G$ is itself soluble. }

{\bf Proposition 3.4.} {\sl Let ${\cal M}$ be an $r$-multiply saturated 
formation and ${\cal N}\subseteq {\cal M}\subseteq {\cal N}^{r+1}$ for some $r\geq 0$.
 Then, for any prime $p$, both formations
 ${\cal M}$ and ${\cal G}_{p}{\cal M}$ are $\Sigma _{r+3}$-closed.}
                                                              
{\bf Proof.}  Let $M$ be 
 the canonical local satellite of $\cal M$.  
 Let ${\cal F}$ be one 
of the formations    
${\cal M}$ or  ${\cal G}_{p}{\cal M}$.    Let $G$ be  any  group such that for 
some subgroups  $H_1,\ldots , H_{r+3}$  of $G$ whose indices 
$|G:H_{1}|, \ldots , |G:H_{r+3}|$  
are pairwise coprime we have   $H_1,\ldots , H_{r+3} \in {\cal F}$. We 
shall prove  $G\in   {\cal F}$.  Suppose that this is false and let 
$G$ be a counterexample with $r + |G|$ minimal.  Let $N$ be a 
minimal normal subgroup of $G$. 

(1)   {\sl $N=G^{\cal F}$ is the only  minimal normal subgroup of $G$
 and  $N\leq O_{q}(G)$ for some prime $q$. Hence if
 ${\cal F}=   {\cal G}_{p}{\cal M}$, then $q\ne p$. }

It is clear that the 
hypothesis holds for $G/N$, so $G/N \in   {\cal F}$ by the choice of $G$. 
Hence $N=G^{\cal F}$ since $G\not  \in   {\cal F}$.  Moreover, 
$N$ is a $q$-group for some prime $q$ since $G$ is soluble by Proposition 3.1.
Finally, if ${\cal F}=   {\cal G}_{p}{\cal M}$ and $p=q$, then 

$$G\in {\cal G}_{p}({\cal G}_{p}{\cal M})={\cal G}_{p}{\cal F}=   {\cal 
F},$$ a contradiction. Hence we have (1).

Since the indices 
$|G:H_{1}|, \ldots , |G:H_{r+3}|$  
are pairwise coprime, in view of (1) we may assume without loss of generality 
that $N\leq H_{i}$ for all $i=2, \ldots , r+3$. 

(2) $C_{G}(N)=N$.

First we show that $N\nleq \Phi (G)$. 
Suppose that $N\leq \Phi (G)$. If $r > 0$, then   $ {\cal F}$  is 
saturated by Lemma 2.8, so  $G \in   {\cal F}$. This contradiction shows 
that $r=0$ and so ${\cal F}=   {\cal G}_{p}{\cal M}$ by Lemma 2.10
 and the choice of $G$.  Hence   $q\ne p$  by (1). Let $O/N=O_{p}(G/N)$ and $P$ be a Sylow $p$-subgroup
 of $O$. Then $G=ON_{G}(P)=NPN_{G}(P)= NN_{G}(P)=N_{G}(P)$ by the Frattini 
Argument since $N\leq \Phi (G)$. Hence in view of (1), $O_{p}(G/N)=1$ and so $G/N\in {\cal M}$ since
 $G/N\in {\cal F}=   {\cal G}_{p}{\cal M}$. But then $G$ is a $p'$-group. 
Hence $H_1,H_2 , H_{3} \in {\cal M}$.  Thus $G \in 
{\cal M}\subseteq  {\cal F}$  by Lemma 2.10.
This contradiction shows that $N\nleq \Phi (G)$. 
But then  $C_{G}(N)=N$ by (1) and \cite[A, Theorem 15.2]{1}.

(3) $r > 0$.

Suppose that $r=0$. Then ${\cal F}=   {\cal G}_{p}{\cal M}$, where   
${\cal M}$ is a formation of nilpotent groups. Since $N\leq H_{2}\in  
{\cal F}$ and, by (2),  $C_{G}(N)=N$,  $O_{p}(H_{2})=1$. Hence $H_{2}$ is 
a $p'$-group. Similarly,  $H_{3}$ is 
a $p'$-group. Hence $G=H_{1}H_{2}$ is 
a $p'$-group. But then $ H_{1}\in  
{\cal M}$, so   $G\in {\cal F}$ by Lemma 2.10.  This contradiction shows 
that we have (3).

(4) {\sl $H_{i}/N\in M(q)$ for all $i=2, \ldots , r+3$. }

Let   $i\in \{2, \ldots  , r +3 \}$.  Then 
  $H_{i}\in {\cal M}$.  Indeed, 
if ${\cal F}=   {\cal G}_{p}{\cal M}$, then $q\ne p$ by (1). On the other 
hand, in view of (2), $C_{G}(N)=N$.
Hence  $O_{p}(H_{i})=1$, which implies that $H_{i}\in {\cal M}$. But then, 
by (2) and Lemma 2.7, $H_{i}/N=H_{i}/C_{H_{i}}(N)\in M(q)$.

(5) $G/N\in M(q)$.

By Lemma 2.11 and \cite[Chapter IV, Proposition 3.8]{1}, $M(q)={\cal G}_{q}{\cal M}_{0}$, where
 ${\cal M}_{0}=l_{r-1}$form  $(G/O_{q', q}(G) | G\in {\cal 
M})$. Since ${\cal M}\subseteq {\cal N}^{r+1}$, $G/O_{q', q}(G) \in 
 {\cal N}^{r}$, so ${\cal M}_{0}\subseteq {\cal N}^{r}$ since   ${\cal M}_{0}$ is an 
$(r-1)$-multiply saturated 
formation. Therefore, the minimality of $r + |G|$ and Claim (4) imply that $G/N\in M(q)$.

{\sl Final contradiction.}  Since $N$ is a $q$-group by (1),  from 
(5) it follows that $G\in {\cal G}_{q}M(q)=M(q)\subseteq  
{\cal M}\subseteq {\cal G}_{p}{\cal M}$. This contradiction completes the 
proof of the  proposition.

{\bf Corollary 3.5} (See \cite[Satz 1.3]{16}).  {\sl  Every saturated formation 
contained in ${\cal N}^{2}$ is $\Sigma _4$-closed.}

{\bf Corollary 3.6.}  
{\sl The class of all soluble groups of nilpotent length at most $r$} ($r\geq 2$) 
{\sl is $\Sigma _{r+2}$-closed.}

{\bf Proof.} It is clear that ${\cal N}^r$ is 
hereditary formation. Moreover, in view of Lemma 2.8,
 ${\cal N}^r$ is an $(r-1)$-multiply saturated 
formation. So ${\cal N}^r$ is $\Sigma _{r+2}$-closed by Proposition 3.4.

{\bf Proof of Theorem A.} 
Suppose that the theorem is false and consider some counterexample $G$
of minimal    order.  
Take a maximal subgroup $M$ of $G$.  
Then by hypothesis all $(n-1)$-maximal subgroups of $M$ are 
$\cal F$-subnormal in $G$, and so they are $\cal F$-subnormal in $M$ 
by Lemmas 2.1(1) and 2.12.  The solubility of $G$ implies 
that  either  $|\pi (M)| =  |\pi (G)|$ or    $|\pi (M)| =  |\pi (G)| - 1$,  
so   $M\in \cal F$  by  the  choice 
of $G$.  
Hence $G$ is an $\cal F$-critical group. 

Since $G$ is soluble, $G$ has a maximal subgroup $T$ with $|G:T|=p^{a}$ for any prime $p$
 dividing $|G|$. On the other hand, $\cal F$ is $\Sigma _{r+3}$-closed by Proposition 3.4.
Hence $|\pi (G)| \leq r+2$.    Moreover, by hypothesis,   $|\pi 
(G)|\geq n+r+1$. Therefore   $n=1$. Thus  all maximal subgroups of  $G$ are $\cal 
F$-normal, so $G/\Phi (G)\in {\cal F}$. But 
 ${\cal F}$ is a saturated formation and hence    $G\in \cal F$. This contradiction
 completes the proof of the result.

{\bf Corollary 3.7} (See \cite[Theorem 6]{13}). {\sl If each $n$-maximal subgroup of a
soluble group $G$
is subnormal, and if $|\pi (G)|\geq n+1$, then $G$ is nilpotent.}

{\bf Corollary 3.8} (See \cite[Theorem A]{19}). {\sl If every
$n$-maximal subgroup of a soluble group $G$
is $\cal U$-subnormal in $G$ and $|\pi (G)|\geq n+2$, then $G$ 
is supersoluble.}

{\bf Corollary 3.9.} {\sl Let $\cal F$ be the class of all groups 
$G$ with $G'\leq F(G)$. If every
$n$-maximal subgroup of a soluble group $G$
is $\cal F$-subnormal in $G$ and $|\pi (G)|\geq n+2$, then $G\in \cal F$.}

{\bf Corollary 3.10.} {\sl If every
$n$-maximal subgroup of a soluble group $G$
is ${\cal N}^r$-subnormal in $G$} ($r\geq 1$) {\sl and $|\pi (G)|\geq n+r$, then 
$G\in {\cal N}^r$.}

\section{Proofs of Theorems B, C, and D}

{\bf Proof of Theorem C.} First suppose that every $2$-maximal subgroup of 
$G$ is $\cal F$-subnormal in $G$. Assume that $G\not \in \cal F$.
We shall show that
$G$ is an $\cal F$-critical group
and $G^{\cal F}$ is a minimal normal subgroup of $G$. 
Let $M$ be a maximal subgroup of $G$
and $T$ a maximal subgroup of $M$.
By hypothesis, $T$ is $\cal F$-subnormal in $G$.
Therefore $T$ 
is $\cal F$-normal in $M$ by Lemma 2.1(1), so $M/T_M\in \cal F$. 
Since $T$ is arbitrary and $\cal F$ is saturated, $M\in \cal F$. 
Consequently, all maximal subgroups of $G$ belong to
$\cal F$.
Hence $G$ is an $\cal F$-critical group.
Then by hypothesis, $G$ is soluble and it  has a normal Sylow $p$-subgroup $G_p\ne 1$
for some prime $p$.
Thus $G_p=G^{\cal F}$   by Lemma 2.5. On the other hand, by Lemma 2.4,  
$G_p/\Phi (G_p)$ is a chief factor of $G$. 

Let $M$ be an $\cal F$-abnormal maximal subgroup of 
$G$. Then $G_p\nleq M$, so
$G=G_{p}M$ and $M=(G_{p}\cap M)G_{p'}=\Phi (G_p)G_{p'}$, where $G_{p'}$
is a Hall $p'$-subgroup of $G$. Assume that $\Phi (G_p)\ne 1$.
It is clear that    $\Phi (G_p)\nleq \Phi (M)$.
Let $T$ be a maximal subgroup of $M$ such that 
$\Phi (G_p)\nleq T$. Then  $M=\Phi (G_p)T$. 
Since $T$ is
$\cal F$-subnormal in $G$,
there is a maximal subgroup  $L$ of $G$ such that 
$T\leq L$ and  $G/L_G\in \cal F$.
Then $G_p\leq L_G$, so  $G=G_{p}M= G_{p}\Phi (G_p)T = G_{p}T\leq L$, a contradiction.
Hence $\Phi (G_p)=1$. Therefore $G_p=G^{\cal F}$
is a minimal normal subgroup of $G$ by Lemma 2.4.

Now suppose that  $G$  is an $\cal F$-critical group
and $G^{\cal F}$ is a minimal normal subgroup of $G$.
Let
$T$ be a $2$-maximal subgroup of $G$ and
$M$ a maximal subgroup of $G$ such that $T$ is a maximal subgroup of $M$.
Since $M\in 
\cal F$, $T$ is $\cal 
F$-subnormal in $M$ by Lemma 2.2. Therefore, if $M$ is $\cal F$-normal in 
$G$, then $T$ is 
$\cal F$-subnormal in $G$ by Lemma 2.1(3). Assume that $M$ 
is $\cal F$-abnormal in $G$. Then $G^{\cal F}\not\leq 
M$. Therefore, since $G^{\cal F}$ is a minimal normal subgroup 
of $G$ by hypothesis, $G=G^{\cal F}\rtimes M$ and  $G^{\cal F}T$ 
is a maximal $\cal F$-normal subgroup of $G$.
 Moreover, since $G$  is an $\cal F$-critical group, $G^{\cal F}T\in {\cal F}$ and hence 
 $T$ is $\cal F$-subnormal in $G^{\cal F}T$
by Lemma 2.2. Hence,  $T$ 
is $\cal F$-subnormal in $G$. The theorem is proved.

From Theorem C and Lemma 2.6 we get

{\bf Corollary 4.1} (See \cite[Theorem 3.1]{19}).
  {\sl Every $2$-maximal subgroup of $G$ is $\cal U$-subnormal in $G$
if and only if  $G$ is an  $\cal U$-critical group
and $G^{\cal U}$ is a minimal normal subgroup of $G$.}

{\bf Proof of Theorem B.} 
First suppose that all $n$-maximal subgroups of $G$ 
are $\cal F$-subnormal in $G$.
We shall show, in this case, that
either $G\in \cal F$ or $G$ is a group of the type II.
Assume that this is false and consider a counterexample $G$
for which $|G|+n$ is minimal. Therefore $A = G^{\cal F}\ne 1$.
Then:

(a) {\sl The hypothesis holds for every maximal subgroup of $G$.}

Let $M$ be a maximal subgroup of $G$.
Then by hypothesis, all $(n-1)$-maximal subgroups of $M$ are 
$\cal F$-subnormal in $G$, and so they are $\cal F$-subnormal in $M$ 
by Lemmas 2.1(1) and 2.12.  Moreover, the solubility of $G$ implies 
that  either  $|\pi (M)| =  |\pi (G)|$ or    $|\pi (M)| =  |\pi (G)| - 1$.

(b) {\sl If $M$ is a maximal subgroup of $G$ and $|\pi (M)|=|\pi (G)|$, then 
$M\in \cal F$.}

In view of hypothesis and Lemmas 2.1(1) and 2.12, all $(n-1)$-maximal subgroups of 
$M$ are 
$\cal F$-subnormal in $M$. Since $|\pi (M)|=|\pi (G)|\geq n+1=n-1 +2$,
$M\in \cal F$ by Theorem A.

(c) {\sl If $W$ is a Hall $q'$-subgroup of $G$
for some $q\in \pi (G)$, then
either $W\in \cal F$ or $W$ is a group of the type II.}

If $W$ is not a maximal subgroup of $G$, 
then there is a maximal subgroup $V$ of $G$ such that $W\leq V$ and
$|\pi (V)|=|\pi (G)|$. By (b), $V\in \cal F$. Hence   $W\in \cal F$ by Lemma 2.12.
Suppose that $W$ is a maximal subgroup of $G$. Then
by (a), the hypothesis holds for $W$, so either $W\in\cal F$ or $W$ 
is a group of the type II by the choice of $G$.

(d) {\sl The hypothesis holds for $G/N$, where $N$  is a minimal normal subgroup of $G$.}
  
If $N$  is not a Sylow subgroup 
of  $G$,  then $|\pi (G/N)| =  |\pi (G)|$.
Moreover,  if $H/N$     is  an  $n$-maximal  subgroup  of  $G/N$,
then  $H$  is  an  $n$-maximal subgroup of $G$.
Therefore $H$  is $\cal F$-subnormal in  $G$.  
Consequently, $H/N$  is $\cal F$-subnormal in  $G/N$ 
by Lemma 2.1(2).  But if $G/N$  has no $n$-maximal subgroups, 
then by the solubility of  $G$, the identity 
subgroup  of  $G/N$  is $\cal F$-subnormal  in  $G/N$
and it  is  the  unique $i$-maximal  subgroup  of  $G/N$  for  some $i < n$ 
with $i < |\pi (G/N)|$.  Finally, 
consider  the  case  that  $N$  is  a  Sylow $p$-subgroup  of  $G$.
Let  $E$ be a Hall $p'$-subgroup of $G$.
It is clear that $|\pi (E)| = |\pi (G)| - 1$ and 
$E$ is a maximal subgroup of $G$.

Let $H/N$ be an  $(n-1)$-maximal subgroup of $G/N$. Then $H$
is an $(n-1)$-maximal subgroup of $G$ and $H=H\cap NE=N(H\cap E)$.
There is a chain of subgroups 
$H=H_{0}<H_{1}<\ldots < H_{n-1}=G$ of 
$G$,
where $H_{i-1}$ is a maximal subgroup of $H_{i}$ ($i=1, \ldots , n-1$).
Then $H_{i-1}\cap E$ is a maximal subgroup of $H_{i}\cap E$, for $i=1, \ldots , n-1$.
Indeed, suppose that for some $i$ there is a subgroup $K$ of $H_{i}\cap E$
such that $H_{i-1}\cap E\leq K\leq H_{i}\cap E$. Then
$(H_{i-1}\cap E)N\leq KN\leq (H_{i}\cap E)N$, so 
$H_{i-1}=H_{i-1}\cap EN\leq KN\leq H_{i}\cap EN=H_{i}$. Whence
 either $KN=H_{i-1}$ or $KN=H_{i}$.
If $KN=H_{i-1}$, then $H_{i-1}\cap E=KN\cap E=K(N\cap E)=K$.
In the second case we have $H_{i}\cap E=KN\cap E=K(N\cap E)=K$.
Therefore $H_{i-1}\cap E$ is a maximal subgroup of $H_{i}\cap E$,
so $H\cap E$ is an $(n-1)$-maximal subgroup of $E$.
Since $E$ is a maximal subgroup of $G$, $H\cap E$ is an $n$-maximal subgroup of 
$G$. Hence $H\cap E$ is $\cal F$-subnormal in $G$ by hypothesis.
Therefore $H/N=(H\cap E)N/N$ is $\cal F$-subnormal in $G/N$ by Lemma 2.1(2).

(e) {\sl $|\pi (G)| > 2$.}

If  $|\pi (G)| = 2$, then $n=1$  and so all maximal subgroups of $G$
are $\cal F$-normal by hypothesis. Hence $G\in \cal F$ since $\cal F$ is a 
saturated formation, a contradiction.

(f) {\sl $G$ is an Ore dispersive group.} 

Suppose that this is false.
Take  a  minimal  normal  subgroup  $N$
of  $G$.  Then by (d),
the  hypothesis holds for  $G/N$,
so either $G/N\in \cal  F$ or $G/N$ is a group of the type II.
Thus, in view of $\cal F\subseteq \cal U$ and the choice of $G$,
$G/N$ is  an  Ore  dispersive  group.   
By Lemma 2.13, the  class  of  all  Ore  dispersive  groups is
a saturated  formation. 
Therefore $N$  is  the  unique  minimal  normal  subgroup  of  $G$
and  $N\nleq \Phi (G)$. Hence $\Phi (G)=1$ and  there is a maximal subgroup 
$L$ of $G$ such that $G = N\rtimes L$ and $L_G = 1$.  
Thus $C_G (N ) = N$ by \cite[A, Theorem 15.2]{1}. 

Since  $G$ is soluble,  $G$ has a normal  maximal 
subgroup  $M$  with       
$|G  : M | = p$  for  some  prime $p$  and
either  
$|\pi (M)| =  |\pi (G)|$ or  $|\pi (M)| =  |\pi (G)| - 1$. 
By (a),  the  hypothesis holds for  $M$.
Therefore, in view of $\cal F\subseteq \cal U$ and the choice of $G$, 
$M$  is an Ore dispersive group.
Denote by $q$ the 
greatest number in $\pi (M)$.
Take a Sylow $q$-subgroup $M_q$
of $M$.
Since $M_q$  is a characteristic subgroup of $M$, $M_q$
is normal in $G$.  
Consider the case $|\pi (M)| = |\pi (G)|$ 
first.   Then  $q$  is  the  greatest  prime  divisor  of  the  order  of  
$G$ and  $M_q\ne 1$. Hence $G/M_q$    is  an  Ore 
dispersive group, and by the maximality of $q$,  so is  $G$.
Suppose now that  $|\pi (M)| =  |\pi (G)| - 1$.
If $q > p$, 
then, as above, we conclude that $G$ is an Ore dispersive group as well.  
Let $p > q$. Then $p$ is the 
greatest prime divisor of  $|G|$.
Since $M_q\ne 1$, it follows that $N\leq M_q$, so
$N$  is a $q$-group.  
In addition, since  $|\pi (G)| > 2$ by (e), 
there is a prime divisor $r$ of the order of $G$ such that $q\ne r\ne p$.
Take a Hall $r'$-subgroup  $W$  of  $G$.
Then $PN\leq W$ for some Sylow $p$-subgroup $P$ of $G$.
Moreover,  by (c), $W$  is an Ore dispersive group. 
Hence $P$  is normal in $W$, and so $P\leq C_G (N)=N$.
The resulting contradiction shows 
that $G$ is an Ore dispersive group. 

(g) {\sl $A$ is a nilpotent group.}

Suppose that this is false.
Let $N$ be a minimal normal subgroup of $G$.
Then by (d),  $(G/N)^{\cal F}=G^{\cal F}N/N\simeq G^{\cal F}/G^{\cal F}\cap N$ 
is  a nilpotent group.   
It  is  known  that  the  class  of  all  nilpotent groups is
a saturated  formation. 
Hence in the case when $G$ has a minimal normal subgroup $R\ne N$
we have $G^{\cal F}/(G^{\cal F}\cap N)\cap (G^{\cal F}\cap R)\simeq G^{\cal F}$
is nilpotent. Thus $N$ is  the  unique  minimal  normal  subgroup  of  
$G$ and $N\leq G^{\cal F}$. If $N\leq \Phi (G)$, then
$G^{\cal F}/G^{\cal F}\cap \Phi (G)\simeq (G^{\cal F}/N)/((G^{\cal F}\cap \Phi (G))/N)$
is nilpotent, so $G^{\cal F}$ is nilpotent by Lemma 2.14.
Therefore  $N\nleq \Phi (G)$. Hence $\Phi (G)=1$ and  there is a maximal subgroup 
$L$ of $G$ such that $G = N\rtimes L$ and $L_G = 1$.  Thus $C_G (N ) = 
N$ by \cite[A, Theorem 15.2]{1}
and $N\ne A$. 

{\sl Case 1: $|\pi (G)| = 3$.}
By hypothesis, either all maximal subgroups of $G$ 
or all its $2$-maximal subgroups are $\cal F$-subnormal in $G$.
In the first case we infer that $G\in \cal F$,
which  contradicts  the choice of $G$.
Hence  all  $2$-maximal  subgroups  of  $G$  are $\cal F$-subnormal.
Since $\cal F\subseteq \cal U$, in view of Lemma 2.6, every $\cal F$-critical 
group has a normal Sylow subgroup.
Whence Theorem C implies that $G$ is an $\cal F$-critical group and $A=G^{\cal F}$
is a minimal normal subgroup of $G$.
Therefore $A=N$, a contradiction.

{\sl Case 2: $|\pi (G)|\geq 4$.}
Assume that 
$N$  is a $p$-group,  and take a Sylow subgroup  $P$
of  $G$ such that $N\leq P$.  Observe that if $N\ne P$, then $L\in \cal F$ by (b),
and so $A=N$, a contradiction.  Hence $N  = P$. 

{\sl Case 2.1: $|\pi (G)|=4$.}

(1) {\sl All $3$-maximal  subgroups  
of $G$ are $\cal F$-subnormal  in  $G$
and $L$ is an $\cal  
F$-critical group.}

Since  $G\not\in \cal F$ and $|\pi (G)|  =  4$,  either  all  $2$-maximal  subgroups  of  $G$  or 
all  its  $3$-maximal  subgroups  are $\cal  F$-subnormal  in  $G$.
In the first case $G$ is an $\cal F$-critical group and 
$A=G^{\cal  F}$ is a minimal normal subgroup of $G$ by Theorem C.
Hence $A=N$, a contradiction. 
Therefore
all $3$-maximal  subgroups  
of $G$ are $\cal  F$-subnormal  in  $G$.
Thus all  second  maximal  subgroups  of $G$ belong to $\cal F$ by Lemma 2.3.
Consequently,   either  $L\in \cal  F$  or $L$ is  an $\cal  F$-critical group. 
But in the first case  $N=A$, a contradiction.
Therefore $L$ is an $\cal  
F$-critical group.

(2) {\sl $L=Q\rtimes (R\rtimes T)$, where 
$Q,R,T$ are Sylow subgroups of $G$, $Q=L^{\cal  F}$ 
is a minimal normal subgroup of $L$, and $G^{\cal F}=PQ$.}

Since $N=P$  is a Sylow $p$-subgroup of $G$ and $|\pi (G)|=4$, $|\pi 
(L)|=3$. Hence in view of (f),  $L=Q\rtimes (R\rtimes T)$, where 
$Q,R,T$ are Sylow subgroups of $G$. Moreover,  $Q=L^{\cal  F}$ by Lemma 
2.5 and  $Q$ is a minimal normal subgroup of $L$ by Theorem C
 since every 2-maximal subgroup of 
$L$ is $\cal F$-subnormal  in  $L$ by (1) and Lemmas 2.1(1) and 2.12.  Finally, since
$G/N\not \in {\cal F}$   and $G/PQ\simeq L/Q\in {\cal F}$,   we have  $G^{\cal F}=PQ$.

(3) {\sl $V= PQR$  is not supersoluble. Hence $V\not\in \cal F$.}

Assume that $V$ is a supersoluble  group.
Since $F (V)$  is a  characteristic subgroup  of $V$
and  $V$   is a  normal  subgroup  of  $G$,
$F(V)$ is normal in $G$. Hence every  Sylow  subgroup  of $F (V)$  is  normal  in $G$.
But $N$  is  the
unique minimal  normal  subgroup  of  $G$. 
Therefore $F (V)= N= P$. Thus $V/P\simeq QR$  is  an  abelian 
group.   
Hence $R$ is normal in $L$ and so
 $R\leq F(L)$. In view of Lemma 2.5, $F(L)=Q\Phi (L)$.
Whence $R\leq \Phi (L)$.
This contradiction shows that 
$V$ is not supersoluble. Thus $V\not\in \cal F$ since
 $ \cal F  \subseteq  \cal U$ by hypothesis.

(4) {\sl $V$ is a maximal subgroup of $G$. Hence $|T|=t$ is a prime.}

If $V$ is not a maximal subgroup of $G$, then
there is a maximal subgroup
$U$ of $G$ such that $V\leq U$ and $|\pi (U)|=|\pi (G)|$.
Hence $U\in \cal F$ by (b), so
$V\in \cal  F$ by Lemma 2.12, a contradiction. 
Therefore 
$V$ is a  normal maximal subgroup of  $G$. Whence
$|T|$ is a prime.

(5) {\sl $|Q|=q$ is a prime and $R=\langle x\rangle$ is a cyclic group.}

Since $V$ is a maximal subgroup of $G$ by (4), 
all $2$-maximal subgroups of $V$ are $\cal F$-subnormal in $V$
by (1) and Lemmas  2.1(1) and 2.12.
Hence, in view of (3), $V$ is an $\cal F$-critical group
by Theorem C.  Therefore, in fact, $V$ is an $\cal U$-critical group by 
(3) since $\cal 
F\subseteq \cal  U$. Hence $QR$ is supersoluble. 
Since $V$ is normal in $G$ and $\Phi (G)=1$, $\Phi (V)=1$.
Therefore $QR$ is a  Schmidt group by Lemma 2.6. Hence $R$ is cyclic
and $Q$ 
is a minimal normal subgroup of $QR$ by Lemma 2.4. 
Whence $|Q|$ is a prime.

(6) {\sl $|R|=r$ is a prime and $C_L(Q)=Q$.}

By (4) and  (5),  $L$  is a supersoluble group. 
Suppose that $|R|=r^b$ is not a prime
and let $M$ be a maximal subgroup of $L$ such that $|L:M|=r$.
Let $W=PM$.
Then $\pi (W)=\pi (G)$, so 
$W\in \cal  F$ by (b) and hence  $W$  is supersoluble. Since $C_{G}(N)=N$, 
$F(W)=P$. Hence $W/P\simeq M$  is abelian.   
It is clear that    $Q\leq M$, so  $M\leq C_{L}(Q)$. Hence 
$T\leq  F(L)$. On the other hand, 
 $F(L)=Q\Phi (L)$ by Lemma 2.5. Therefore $T\nleq  F(L)$. This 
contradiction shows that  $|R|=r$  and   so
$C_L(Q)=Q$ by Lemma 2.5.

(7) {\sl $1\ne C_G(x)\cap PQ=P_1\leq P$.}

Suppose that
$C_G(x)\cap PQ=1$. Then by the Thompson's theorem \cite[Theorem 10.5.4]{18}, $PQ$ is a nilpotent group,
so $Q\leq C_G(P)=P$,
a contradiction. Thus $C_G(x)\cap PQ\ne 1$.
Suppose that $q$ divides $|C_G(x)\cap PQ|$.
Then, by (5),  for some $a\in P$ we have $Q^a\leq C_G(x)\cap PQ$,
so $\langle Q^a, RT\rangle \leq N_G(R)$.
Hence if $E$ is a Hall $p'$-subgroup of $N_G(R)$, then $E\simeq L$.
Therefore $L$ has a normal $r$-subgroup, so $C_L(Q)\ne Q$, a contradiction.
Thus $C_G(x)\cap PQ=P_1\leq P$.

{\sl Final contradiction for Case 2.1.}
Let $D=\langle P_1, RT\rangle$. Then $D\leq N_G(R)$. If $q$ divides $|D|$, 
then, as above, we have $C_L(Q)\ne Q$. Thus $D\cap Q^a=1$ for all $a\in P$.
Moreover, if $P\leq D$, then
$PR=P\times R$ and $R\leq C_G(P)=P$. Therefore $P\nleq D$ and
$D$ is not a maximal subgroup of $G$. Hence  $D$ 
is a $k$-maximal subgroup of $G$ for some $k\geq 2$.
Then
there is a $3$-maximal subgroup $S$ of $G$ such that $RT\leq S\leq D$.
By hypothesis, $S$ is $\cal  F$-subnormal in $G$. Hence  at 
least one of the maximal subgroups $L$ or $PRT$ is $\cal F$-normal in $G$, contrary to (2). 

{\sl Case 2.2: $|\pi (G)| > 4$.}
If $\pi (L) = \{p_1, \ldots , p_t\}$, then $t > 3$.
Let $E_i$ be a Hall $p_i'$-subgroup of $L$ and 
 $X_i=PE_i$. We shall show that  $E_i\in \cal F$ for all $i=1, \ldots , t$. 
By (c), either $X_i\in \cal  F$ or 
 $X_i$ is a group of the type II, for $i=1, \ldots , t$. In the former 
case we have $E_{i} \simeq X_{i}/P\in \cal  F$.  Assume that 
 $X_i$ be a group of the type II.
Then $X_i^{\cal F}$ is nilpotent, so 
$X_i^{\cal  F}\leq F(X_i)$.
But since $P$ is normal in $X_i$ and $C_G(P)=P$, 
$F(X_i)=P$. 
Hence $X_i^{\cal F}=P$, so $E_i\in \cal F$.
Since $t > 3$, Proposition 3.4 implies that then $L\in \cal  F$.
Therefore $A=N$, a contradiction.  Hence we have (g).

(h) {\sl $A$ is a Hall subgroup of $G$.}

Suppose that this is false.
Since $G$ is Ore dispersive by (f), for the greatest prime divisor $p$ 
of $|G|$ the Sylow $p$-subgroup $P$ is normal in $G$. 
Assume that $P$ is not a minimal normal subgroup of $G$.
Then there is a maximal subgroup $M$ of $G$ such that
$G=PM$ and $P\cap M\ne 1$.
Since $|\pi (M)|=|\pi (G)|$, $M\in \cal  F$ by (b). 
Hence $G/P\simeq M/M\cap P\in \cal F$, so $A=G^{\cal  F}\leq P$.
Suppose that $\Phi (P)\ne 1$. Let $N$ be a minimal normal subgroup of $G$ 
such that $N\leq \Phi (P)$. 
By (d), the hypothesis holds for $G/N$,
so either $G/N\in \cal F$ or $G/N$ is a group of the type II by the 
choice of $G$.
If $G/N\in \cal F$, then $A=N\leq \Phi (P)$. Since $P$ is normal in 
$G$, $\Phi (P)\leq \Phi (G)$. Thus $A\leq \Phi (G)$ and so $G\in 
\cal F$, a contradiction. Hence $G/N$ is a group of the type II.
Therefore $AN/N=G^{\cal  F}N/N=(G/N)^{\cal F}$ 
is  a  Hall  subgroup  of  $G/N$.
Consequently, $AN=P$. Hence $A\Phi (P)=P$, so $A=P$, a contradiction.
Thus $\Phi (P)=1$.
By Maschke's 
theorem,  $P=N_1\times \ldots \times N_k$  is  the  direct  product  of  
some minimal  normal  subgroups  of  $G$.
If $N_1\ne P$,  then 
$G/N_1\in \cal F$  and $G/N_2\in \cal F$ by Theorem A.
Consequently, so is $G$.
This contradiction shows that
$P$ is a minimal normal subgroup of $G$. 

By (d), the hypothesis holds for $G/P$,
so either $G/P\in \cal  F$ or $G/P$ is a group of the type II by the 
choice of $G$.
If $G/P\in \cal  F$, then $A=P$, a contradiction.
Hence $G/P$ is a group of the type II.
Therefore $AP/P=G^{\cal  F}P/P=(G/P)^{\cal  F}$ 
is  a  Hall  subgroup  of  $G/P$.
If $P\leq A$, then $A=P\rtimes A_{p'}$, where $A_{p'}$
is a Hall $p'$-subgroup of $A$. But since $A_{p'}\simeq A/P$
and $A/P$ 
is  a  Hall  subgroup  of  $G/P$, $A$ is a Hall subgroup of $G$.
Therefore $P\cap A=1$, so $A$ is a Hall subgroup of $G$ since $AP/P\simeq A/A\cap P\simeq A$.

(i) {\sl $A$ is either of the form $N_1\times \ldots \times N_t$, 
where each $N_i$ is a minimal normal subgroup of $G$, which is a Sylow 
subgroup of $G$,  for $i=1, \ldots , t$, or a Sylow $p$-subgroup of $G$ of exponent $p$
for some prime $p$ and the commutator subgroup, the Frattini subgroup, and 
the center of $A$ coincide, while $A/\Phi (A)$ is an $\cal F$-eccentric chief factor
of $G$.}

Suppose that $A$ is not a minimal normal subgroup of $G$.
Take a Sylow $p$-subgroup $P$  of $A$,
where $p$  divides  $|A|$.
Claims (g) and (h) imply that $P$  is a normal 
Sylow  subgroup  of  $G$.
Let  $N$  be  a  minimal  normal  
subgroup  of  $G$  with  $N\leq P$.
First suppose 
that $N\leq \Phi (G)$, and take a maximal subgroup $M$
of $G$ with $P\nleq M$.  Then $M\in \cal F$ by 
(b). Therefore $G/P\simeq M/M\cap P\in \cal  F$.
In this case $A = P$. Moreover, if $S$ is a maximal subgroup
of $G$ such that $P\nleq S$, then $S\in\cal  F$.
Observe also that for every maximal subgroup $X$ of $G$
with $P\leq X$ we have $X$ is $\cal F$-subnormal in $G$. 
Thus, by Lemma 2.4, $A=G^{\cal  F}$ satisfies condition II(1). 

Suppose  that for  every  minimal  normal  subgroup $R$ of  $G$
such that $R\leq A$  we have $R\nleq \Phi (G)$.
Then there is a maximal subgroup $L$ of $G$ such that $G=N\rtimes L$. 
If $N\ne P$,  then  $L\in \cal  F$
by  (b).  Therefore $A=N$,  a contradiction. 
Consequently, all Sylow subgroups of $A$
are minimal normal subgroups of $G$.
Therefore $A=N_1\times \ldots \times N_t$, where $N_i$
is a minimal normal subgroup of $G$, for $i = 1,\ldots , t$.

(j) {\sl  Every 
$n$-maximal subgroup   of $G$
belongs to $\cal F$ and induces on the Sylow $p$-subgroup of 
$A$ the  automorphism group which is contained in $F(p)$
 for every prime divisor $p$ of $|A|$.}

Let $H$ be any $n$-maximal subgroup   of $G$. 
Suppose that $H$ is a maximal subgroup of $V$, where $V$ is an 
$(n-1)$-maximal subgroup of $G$. Since $V\in 
\cal  F$ by Lemmas 2.3 and 2.12, $H\in \cal  F$.

Let $E=AH$.  Since $A$ 
is normal in $E$ and $A$ is nilpotent by (g),
$A\leq F(E)$. Whence $E=F(E)H$. Since $H$  is
$\cal  F$-subnormal in $G$, $H$ is $\cal F$-subnormal in $E$  
by Lemmas 2.1(1) and 2.12. Moreover, $H \in \cal  
F$. Therefore $E\in \cal F$ by Lemma 2.15.
Let $P$ be a Sylow $p$-subgroup of $A$
and $K/L$ a chief factor of $E$ such that $1\leq L<K\leq P$.
Since $E\in \cal  F$, $E/C_E(K/L)\in F(p)$. Hence $P\leq 
Z_{\cal F}(E)$, so $E/C_E(P)\in F(p)$ by Lemma 2.7. Then $H/C_H(P)=H/C_E(P)\cap H
\simeq HC_E(P)/C_E(P)\in F(p)$.

Now suppose that either  $G\in \cal F$ or
$G$ is a group of type II.
If $G\in \cal  F$, 
then every subgroup of $G$ is $\cal  F$-subnormal in $G$ by Lemma 2.2.
Let $G$ be a group of type II. 
Take an $n$-maximal subgroup $H$ of $G$. Put $E=G^{\cal F}H$. 
Let $P$ be a Sylow $p$-subgroup of $G^{\cal  F}$
and $K/L$ a chief factor of $E$ such that $1\leq L<K\leq P$.
By hypothesis, $H/C_H(P)\in F(p)$, so
$H/C_H(K/L)\simeq (H/C_H(P))/(C_H(K/L)/C_H(P))\in F(p)$.
Since $G^{\cal F}$ is normal in $E$ and $G^{\cal  F}$
is nilpotent,  $G^{\cal F}\leq F(E)\leq 
C_E(K/L)$.
Hence $$E/C_E(K/L)=E/C_E(K/L)\cap E=E/C_E(K/L)\cap G^{\cal  F}H=
E/G^{\cal  F}(C_E(K/L)\cap H)=$$ $$=G^{\cal  F}H/G^{\cal  F}C_H(K/L)\simeq
H/G^{\cal  F}C_H(K/L)\cap H=H/C_H(K/L)(G^{\cal  F}\cap 
H)=H/C_H(K/L),$$ so $E/C_E(K/L)\in F(p)$ since $F(p)$ is hereditary by Lemma 2.12
 and \cite[Proposition 3.1.40]{15}.
Then $P\leq Z_{\cal  F}(E)$, whence $G^{\cal  F}\leq Z_{\cal F}(E)$.
Thus $E/Z_{\cal  F}(G)\in \cal  F$.  Hence $E\in 
\cal  F$, so $H$ is an
$\cal  F$-subnormal subgroup of $G^{\cal F}H=E$. Since 
$G^{\cal F}\leq G^{\cal F}H$,  
$G^{\cal F}H$ is $\cal  F$-subnormal in $G$ by Lemma 2.1(4). 
Consequently, in view of Lemma 2.1(3), $H$ is $\cal F$-subnormal in $G$. 
The theorem is proved.

{\bf Corollary 4.2} (See \cite[Theorem B]{19}). 
{\sl Given a soluble group $G$ with $|\pi (G)|\geq n+1$, all
$n$-maximal subgroups of $G$ 
are $\cal  U$-subnormal in $G$ if and only if $G$ is a group of one of 
the following types:}

I. {\sl $G$ is supersoluble.}

II. {\sl $G=A\rtimes B$, where $A=G^{\cal U}$  and  $B$ are Hall subgroups of $G$, 
while $G$ is Ore dispersive and satisfies the following: }

(1) {\sl $A$ is either of the form $N_1\times \ldots \times N_t$, 
where each $N_i$ is a minimal normal subgroup of $G$, which is a Sylow 
subgroup of $G$,  for $i=1, \ldots , t$, or a Sylow $p$-subgroup of $G$ of exponent $p$
for some prime $p$ and the commutator subgroup, the Frattini subgroup, and 
the center of $A$ coincide, every  chief factor of $G$ below $\Phi (G)$
is cyclic, while $A/\Phi (A)$ is a noncyclic chief factor
of $G$;  }

(2) {\sl for every prime divisor $p$ of the order of $A$ 
every $n$-maximal subgroup  $H$ of $G$ is supersoluble and induces on the Sylow $p$-subgroup of 
$A$ an automorphism group which is 
an extension of some $p$-group by abelian group of exponent dividing $p-1$.}

{\bf Proof of Theorem D.}
Assume that this is false and consider a counterexample $G$
for which $|G|+n$ is minimal. 

(a) {\sl $G$ has a unique minimal normal subgroup $N$ such that
$C_G(N) = N$ and $N$ is not a Sylow subgroup of $G$.}

Let $N$ be a minimal normal subgroup of $G$.
Then the  hypothesis holds  for  $G/N$ (see 
Claim (d) in the proof of Theorem B). 
Consequently,  $G/N$ is $\phi$-dispersive for some ordering $\phi$ 
of $\mathbb P$ by the choice of $G$.  
Therefore $N$
is not a Sylow subgroup of $G$.
Moreover, by Lemma 2.13, 
$N\nleq \Phi (G)$.  Therefore 
$G$ has  a  maximal  subgroup  $M$
such  that  $G=N\rtimes M$. By Lemmas 2.1(1) and 2.12
all $(n - 1)$-maximal subgroups of $M$  are $\cal F$-subnormal in 
$M$. Moreover,  $|\pi (M)| = |\pi (G)|$.
Therefore Theorem B 
implies that $G/N\simeq M$  is an Ore dispersive group.
Hence in the case when $G$ has a minimal normal subgroup $R\ne N$
we have $G/N\cap R\simeq G$
is an Ore dispersive group.
Thus $N$ is  the  unique  minimal  normal  subgroup  of  
$G$, and so  $C_G(N) = N$ by \cite[A, Theorem 15.2]{1}. 

(b) {\sl If $W$ is a Hall $q'$-subgroup of $G$ for some $q\in \pi (G)$, then
$W$ is $\phi$-dispersive for some ordering $\phi$ 
of $\mathbb P$.}

If $W$ is not a maximal subgroup of $G$, 
then there is a maximal subgroup $V$ of $G$ such that $W\leq V$ and
$|\pi (W)|=|\pi (G)|$. By hypothesis, every $(n-1)$-maximal subgroup of $V$
is $\cal F$-subnormal in $G$, so it is $\cal F$-subnormal in $V$ 
by Lemmas 2.1(1) and 2.12. Then, in view of Theorem B, $V$ is Ore dispersive.
Hence $W$ is Ore dispersive.
Suppose that $W$ is a maximal subgroup of $G$. Then
$|\pi (W)|=|\pi (G)|-1$ and every $(n-1)$-maximal subgroup of $W$
is $\cal F$-subnormal in $W$ in view of hypothesis and Lemmas 2.1(1) and 2.12.  
Therefore $W$ is $\phi$-dispersive
for some ordering $\phi$ of $\mathbb P$ by the choice of $G$.

(c) {\sl $|\pi (G)| > 2$.}

Suppose  that  $|\pi (G)| =  2$.
Then  by  hypothesis,  either  all  maximal  
subgroups of  $G$ or all its $2$-maximal subgroups 
are $\cal  F$-subnormal in  $G$.
Therefore every maximal subgroup 
of  $G$  belongs to $\cal  F$ in view of Lemmas 2.3 and 2.12.
Consequently,  either  $G\in \cal  F$  or $G$ is  an
$\cal  F$-critical group.
Since $\cal  F\subseteq \cal U$, $G$ is either a supersoluble 
group or an $\cal U$-critical group.
Therefore, in view of Lemma 2.6, $G$ is $\phi$-dispersive for some ordering $\phi$
of $\mathbb P$, a contradiction.

{\sl Final contradiction.} Suppose  that  $N$ is  a $p$-group,  and  
take  a  prime  divisor  $q$  of  $|G|$  such that $q\ne p$.
Take 
a Hall $q'$-subgroup  $E$  of  $G$.
Then $N\leq E$.
By (b), $E$ is $\phi$-dispersive
for some ordering $\phi$ of $\mathbb P$.
Consequently, 
some Sylow subgroup $R$ of $E$  is normal in $E$.
Furthermore, if $N\nleq R$, then  $R\leq C_G (N) = N$. 
Hence  $R$  is  a  Sylow  $p$-subgroup  of  $E$.
It  is  clear  also  that  $R$  is  a  Sylow  $p$-subgroup  of  $G$  and  
$(|G: N_G (R)|, r)  =  1$  for  every  prime  $r\ne q$.
Since  $|\pi (G)|  >  2$ by (c),
$R$  is 
normal in $G$. 
Hence $G$ is $\phi$-dispersive for some ordering $\phi$ 
of $\mathbb P$, a contradiction.
The theorem is proved.

{\bf Corollary 4.3} (See \cite[Theorem C]{19}). {\sl 
If every $n$-maximal subgroup of a soluble group $G$ is
$\cal  U$-subnormal in $G$ and $|\pi (G)|\geq n$, then $G$ is 
$\phi$-dispersive for some ordering $\phi$  of the set of all primes. }

Finally, note that there are examples which show
that the restrictions  on  $|\pi (G)|$
in  Theorems  A,  B,  and  D  cannot  be  weakened.

\end{document}